\documentstyle[12pt,epsf,amscd]{amsart}
\newtheorem{Lemma}{Lemma} [section]
\newtheorem{Proposition}{Proposition} [section]
\newtheorem{Theorem}{Theorem} [section]
\newtheorem{Corollary}{Corollary} [section]
\newtheorem{Definition}{Definition} [section]
\newtheorem{Notation}{Notation} [section]
\newtheorem{Example}{Example} [section]
\newtheorem{Remark}{Remark} [section]
\def\proof{\par{\it Proof}. \ignorespaces}
\def\endproof{{\ \vbox{\hrule\hbox{%
   \vrule height1.3ex\hskip0.8ex\vrule}\hrule }}\par}
\newenvironment{Proof}{\proof}{\endproof}

\begin{document}

\title{Twisted Tomei manifolds and Toda lattice$^*$
}

\author{Luis Casian$^a$ \ \
 and  \ Yuji Kodama$^b$ \\ }

\thanks{$a$) Department of Mathematics, The Ohio State University,
Columbus, OH 43210\endgraf
{\it E-mail address\/}: casian@@math.ohio-state.edu
\endgraf $b$) Department of Mathematics, The Ohio State University,
Columbus, OH 43210\endgraf
{\it E-mail address\/}: kodama@@math.ohio-state.edu
\endgraf $^*$Supported by NSF grant DMS0071523.}

\begin{abstract} 
This paper begins with an observation that the  isospectral leaves of the 
signed  Toda lattice as well as the Toda flow itself may  be constructed from 
the Tomei manifolds by cutting and pasting along certain chamber  walls inside 
a polytope.  It is also observed through examples that although  there is some 
freedom in this procedure of cutting and pasting the manifold and the flow,  
the choices  that  can be made are not arbitrary.   We proceed to describe a 
procedure that begins with an action of the Weyl group on a set of signs; it 
 uses the  Convexity Theorem in \cite{BFR:90} and combines the resulting 
polytope with the chosen Weyl group action to paste together  a compact 
manifold. This manifold which is obtained carries an action of the Weyl group and a Toda lattice
 flow which is related to this action. This construction gives rise to a large 
family of compact manifolds which is parametrized by twisted sign actions of
 the Weyl group. For example, the trivial action gives rise to
   Tomei manifolds and  the standard action of the Weyl group on the connected
    components of a split Cartan subgroup of a split semisimple real Lie group gives rise to the isospectral leaves of the signed Toda lattice. This clarifies
 the connection between the polytope in the Convexity Theorem  and the topology 
of the compact smooth manifolds arising from the  isospectral leaves of a Toda 
flow.  Furthermore, this allows  us to give a uniform treatment to  two very 
different cases that have been studied extensively in the literature producing new 
cases to look at.  Finally we describe the unstable manifolds of the Toda 
flow for these  more general manifolds and determine which of these give rise
 to cycles. 
 
\end{abstract}

\maketitle

\markboth{LUIS CASIAN AND YUJI KODAMA}
  {TWISTED TOMEI MANIFOLDS}

\section{introduction }
\renewcommand{\theequation}{1.\arabic{equation}}\setcounter{equation}{0}
\renewcommand{\thefigure}{1.\arabic{figure}}\setcounter{figure}{0}

In \cite{tomei:84}, Tomei constructed a compact and orientable smooth
manifold as an iso-spectral real manifold generated by the Toda lattice 
equation on the set of the tridiagonal symmetric matrices. Let us begin with
a brief description of this manifold: Let $Z$ be the set of $(l+1)\times (l+1)$ tridiagonal trace zero matrices,
\begin{eqnarray}
\label{X}
Z = \left\{ \left(
\begin{array}{lllll}
a_1 & b_1 & 0 & \cdots & 0 \\
b_1 & a_2 & b_2 & \cdots & 0 \\
\vdots & \ddots & \ddots & \ddots & \vdots \\
0 & \cdots & \ddots & a_{l} & b_{l} \\
0 & \cdots & \cdots & b_{l} & a_{l+1} \\
\end{array}
\right) \quad : ~~
\sum_{i=1}^{l+1} a_i=0, ~ b_i\ne 0 ~~ \right\}.
\end{eqnarray}
According to the signs of $b_i$'s, the set $Z$ is a disjoint union of
$2^{l}$ connected components;
\begin{equation}
\label{Z}
Z=\displaystyle{\bigcup_{\epsilon\in {\cal E}} Z_{\epsilon}},
\end{equation}
where the set of signs is defined by
\begin{equation}
\label{sign}
{\cal E}:=\{~(\epsilon_1,\cdots,\epsilon_l) \in \{\pm 1\}^l~:~\epsilon_k={\rm sign}(b_k), ~ k=1,\cdots,l~\}.
\end{equation}
An isospectral leaf is a subset of $Z$ consisting of tridiagonal marices
with fixed eigenvalues $\lambda_1,\cdots,\lambda_{l+1}$ with
$\sum_{k=1}^{l+1}\lambda_k=0$.  We consider in this paper only the case with
distinct eigenvalues. Any point of
the leaf is obtained from the diagonal matrix 
$\Lambda={\rm diag}(\lambda_1,\cdots,\lambda_l)$ by conjugation, that is,
they lies on an orbit of the adjoint action of the group of orthogonal
 matrices on the set of symmeric matrices,
\begin{equation}
\label{orbit}
{\cal O}_{\Lambda}=\{~{\rm Ad}_p\Lambda~:~p\in SO(l+1)~\}.
\end{equation}
The Toda lattice 
is given by the following matrix equation for $X(t)\in Z$ with a parameter
$t\in {\Bbb R}$,
\begin{equation}
\label{toda}
{dX\over dt}=[P, X]
\end{equation}
where the matrix $P$ is the skew symmetrizaion of $X$, i.e.
$P=X_+-X_-$ with $X_+ (X_-)$ the upper (lower) triangular part of $X$.
The solution defines an orbit in ${\cal O}_{\Lambda}$,
\begin{equation}
\label{solution}
X(t)={\rm Ad}_{Q(t)} X(0),
\end{equation}
where the matrix function $Q(t)\in SO(l+1)$ is given by the QR-factorization 
of the matrix $\exp(tX(0))$,
\begin{equation}
\label{QR}
\exp(tX(0))=Q(t)R(t),
\end{equation}
with an upper triangular matrix $R(t)$.
The isospectral leaf can be also expressed as the inverse image of
the symmetric polynomials $\gamma_1,\cdots,\gamma_l$ of the eigenvalues,
\begin{equation}
\label{zgamma}
Z(\gamma)={\cal I}^{-1}(\gamma)\cap Z
\end{equation}
where $\gamma=(\gamma_1,\cdots,\gamma_l)\in {\Bbb R}^l$ is
the image of the Chevalley invariants ${\cal I}=(I_1,\cdots,I_l):
Z\rightarrow {\Bbb R}^l$ defined
 in the characteristic polynomial of the matrix $X\in Z$,
\begin{equation}
\nonumber
{\rm det}(\lambda I-X)=\lambda^{l+1}-\sum_{k=1}^l (-1)^kI_k\lambda^{l-k}
=0.
\end{equation}
It is then an easy exercise that the $I_k$'s are the polynomial
functions of $a_i$'s and $b_j^2$'s, i.e. $I_k(a_i, b_j)=I_k(a_i,-b_j)$
for $k=1,\cdots,l$.
This implies that each connected component $Z_{\epsilon}(\gamma)$
labeled by the signs of $b_i$'s is just a copy of one single component,
say $Z^+(\gamma)$ having all $b_i>0$ (in Section 2, the topology
of $Z^+(\gamma)$ and its closure ${\overline Z}^+(\gamma)$ will be discussed in detail).
Then Tomei's manifold is considered as a smooth compactification of those 
components under a prescribed gluing through the boundaries, in fact,
the manifold has a $CW$-decomposition with the cells marked by 
sequences of signs and zeros of $b_i$'s.
The smoothness of the manifold was shown by giving local charts using the Toda flow (Lemma 2.1, 2.2 in \cite{tomei:84}). This can be also shown by considering
the manifold as a Morse complex associated with the Toda flow
which gives a {\it gradient-like}
 flow on the manifold. The Morse function 
 $f:{\overline Z}(\gamma)\rightarrow {\Bbb R}$ is given by 
\begin{equation}
\label{morse}
f(X)=\displaystyle{\sum_{i=1}^{l}(l-i+1)a_i},
\end{equation}
where ${\overline Z}(\gamma)$ is the closure of the
isospectral set
(see also \cite{bloch:92} for the Morse functions for the 
generalized Toda lattices on the semisimple Lie algebras).
All the critical points of the function are then given by 
the diagonal matrices with the $(l+1)!$ permutation of
$(\lambda_1,\cdots,\lambda_{l+1})$, i.e. all $b_i=0$.
One can also show that near the critical point $(\lambda_{\sigma(1)},
\cdots,\lambda_{\sigma(l+1)})$ with a permutation $\sigma\in S_{l+1}$,
the symmetry group of order $l+1$, the Morse function takes the form,
\begin{equation}
\nonumber
f(X)=\displaystyle{\sum_{i=1}^l (l-i+1)\lambda_{\sigma(i)}
-\sum_{i=1}^{l} {1\over \lambda_{\sigma(i)}-\lambda_{\sigma(i+1)}}
b_i^2 + o(b_k^2)}.
\end{equation}
This formula is useful for identifying all the handle bodies in
the Morse decomposition of the manifold. However, our primary interest in
this paper is not Morse theory; we will leave that for a future 
communication.

In this paper, we consider an extension of this Tomei manifold
based on a twisted action of the Weyl group, the symmetry group
for the example of this section.  We start with the theorem due to Bloch et al  
in the case of a Lie algebra of type $A_l$ \cite{BFR:90} 
(which is summarized in its general Lie theoretic setting in subsection \ref{moment} below):

\begin{Remark} The imbedding $\iota$.
\label{RemBFR}
\end{Remark}
The isospectral manifold ${\overline Z}^+(\gamma)$ in (\ref{zgamma})
is the closure of a generic noncompact torus orbit in ${\cal O}_{\Lambda}$
of (\ref{orbit}), with an action which does not agree with usual action of a split Cartan
subgroup (the real diagonal matrices of determinant one) on  ${\overline Z}^+(\gamma)$. This poses a difficulty in the application of a Convexity Theorem of Atiyah in \cite{atiyah:82} to this situation. However there is another imbedding $\iota$
described in \cite{BFR:90}  which corrects this difficulty.

The following is an $A_l$ version 
of the main theorem in  \cite{BFR:90}:

\begin{Theorem}
\label{BFRconvexity}

Theere is an imbedding 
$\iota: {\overline Z}^+(\gamma)\hookrightarrow {\cal O}_{\Lambda}$
such that its image is the closure of an orbit under the usual action of  real diagonal matrices of determinant one (a split Cartan subgroup).
Using Atiyah's theorem
\cite{atiyah:82}, the image of the moment map 
$\Gamma:=J\circ \iota ({\overline Z}^+(\gamma) )$  where
$J: X\in{\cal O}_{\Lambda} \rightarrow {\rm diag}X$ is orthogonal projection,
 is the convex polytope with vertices at the critical
points of the Morse function (\ref{morse}), that is, of the
orbit of the Weyl group $S_{l+1}$ of the Lie algebra $\frak{sl}(l+1)$.
\end{Theorem}
The theorem can be applied to each isospectral leaf ${\overline Z}_{\epsilon}
(\gamma)$ with different signs $\epsilon\in {\cal E}$ of (\ref{sign}).
Then Tomei's manifold $M_T$ is the union of those convex polytopes 
attached along
the boundaries, and the Weyl group acts on each polytope independently.
Thus the Tomei manifold can be described as the union of
the polytopes $\Gamma_{\epsilon}$ associated with the isospectral
leaves ${\overline Z}_{\epsilon}$,
\begin{equation}
\label{tomei}
M_T=\displaystyle{\bigcup_{\epsilon\in {\cal E}} \Gamma_{\epsilon}}.
\end{equation}
The main idea of the present paper is to consider a twisted 
action of the Weyl group on the signs in ${\cal E}$, and based on this group action
with specific sign change, we construct a smooth gluing between the 
Weyl chambers inside the polytopes with different signs. This is consistent 
with the Toda flow, since the Morse function generating the flow depends
only on $\{b_i^2 : i=1,\cdots,l\}$ and so does the Toda flow.
Note that a point of the change of signs occurs at the same
point on the same chamber walls of the polytopes with corresponding signs, 
so that the resulting flow stays smooth and provides a smooth gluing
along the chamber walls.

\begin{figure}
\epsfysize=9cm
\centerline{\epsffile{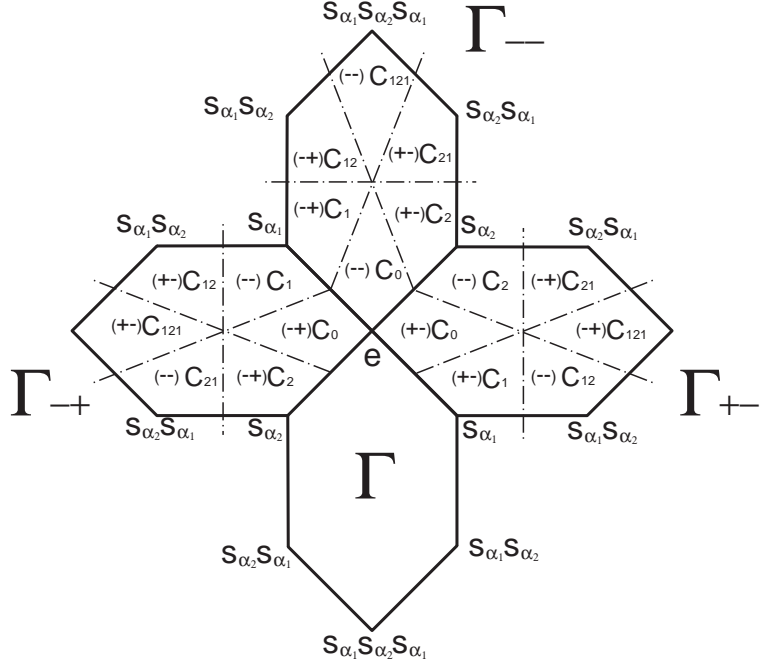}}
\caption{The polytopes $\Gamma_{\epsilon}$ with a standard
$S_3$-action. Chambers $C_{w^{-1}}$ are marked as $C_0$ for $w=e$, $C_1$ 
for $w={s_{\alpha_1}}$ and so on.}
\label{fig:0}
\end{figure}

Let us explain the result with an explicit example of $l=2$.
The polytope is then the hexagon with the six vertices
corresponding to the permutation of $\lambda_1,\lambda_2,\lambda_3$,
which are the critical points of the Morse function (\ref{morse}).
We here assume $\lambda_1>\lambda_2>\lambda_3$.
Each Weyl chamber contains one critical point which is marked by
the corresponding element of the Weyl group. In particular,
the critical point in the dominant chamber is associated with
the fixed point $X_0={\rm diag}(\lambda_1,\lambda_2,\lambda_3)$ of the Toda flow.
The polytope $\Gamma$ is now devided into the six Weyl chambers
$C_w, w\in S_3$, and we let denote $\Gamma_{\epsilon}$
the polytope which has the sign $\epsilon$ in the dominant
chamber and has a prescribed gluing along the inner walls of the chambers,
\begin{equation}
\label{gammaep}
\Gamma_{\epsilon}=\displaystyle{\bigcup_{w\in S_3}\{w\epsilon\}\times
{\overline C}_{w^{-1}}}.
\end{equation}
where ${\overline C}_w$ is the closure of $C_w$ (see Figure \ref{fig:0}).
Here $C_{w}$ describes the chamber containing the critical
point $wx_0$ with $x_0=J(X_0)$.
The prescribed gluing is given by the $S_3$-action on 
${\cal E}=\{\pm 1\}^2$, and
our smooth compact manifold has an expression,
\begin{equation}
\label{manifold}
 M(\delta)=\bigcup_{\epsilon \in \{\pm 1\}^2}\Gamma_{\epsilon},
 \end{equation}
 where $\delta$ denotes the Dynkin diagram $\bullet -\bullet$ having
 a sign $+$ or $-$ on the edge, which specifies the $S_3$-action on the sign.
 This marked Dynkin diagram is defined in Definition 3.3.
If the action is trivial, $w\epsilon=\epsilon, \forall w\in S_3, 
\forall \epsilon
\in {\cal E}$ (the corresponding marked Dynkin diagram is
$\delta=\bullet \overset{+}{-}\bullet$), then we have the Tomei manifold 
(\ref{tomei})
which is topologically equivalent to the connected sum
of two torus.
If the action is a standard one with $s_{\alpha_i}\in S_3:\epsilon_j
\mapsto \epsilon_j\epsilon_i^{C_{j,i}}$ with $(C_{i,j})$ the Cartan
matrix for $A_2$ (then $\delta=\bullet\overset{-}{-}\bullet$), we obtain the compact manifold corresponding
to the indefinite Toda equation discussed in \cite{kodama:98, casian:99}
which is equivalent to the connected sum of two Klein bottles.

The paper is organized as follows:

In Section 2, we present the background information on the Toda
lattice on the variety of ad-diagonalizable Jacobi elements with
nonnegative signs.  Then the convex polytope of the Weyl orbit
associated with the
Toda lattice is introduced based on the Convexity Theorem due to Bloch et al
\cite{BFR:90}. We also present the cell decomposition consisting
of the Weyl chambers and their walls.

In Section 3, we define a twisted action of the Weyl group on the set of signs
${\cal E}=\{\pm 1\}^l$. We also introduce a marked Dynkin diagram
where all the edges are assigned the signs which parametrizes a 
twisted sign action of the Weyl group (Definition 3.3). Then we obtain the
number of compatible twisted action of the Weyl group (Proposition 3.1).

In Section 4, we construct a twisted Tomei manifold associated with a
marked Dynkin diagram by defining gluing maps along both external
and internal Weyl chamber walls.
The manifold is then shown to have the structure of a smooth and compact
manifold with a Weyl group action (Proposition 4.1). 
A cell decomposition of the manifold is also obtained in terms of 
signed colored Dynkin diagrams introduced in \cite{casian:99}.
Then we conclude that the manifold is not orientable if the marked
Dynkin diagram has at least one negative sign on the edge
(Proposition 4.3), and the homology groups over ${\Bbb Z}/2{\Bbb Z}$
of all the twisted Tomei manifolds are the same as that of the original
Tomei manifold (Theorem 4.1).

We then expect that the universal cover of
the twisted Tomei manifold $M(\delta)$ is diffeomorphic to
${\Bbb R}^l$, hence $M(\delta)$ is also aspherical.
Also we would like to mention that
the homology group $H_*(M(\delta), {\Bbb Z})$ can be computed by the cell decomposition given in section 4.2, but an explicit computation 
remains open. Since the twisted Tomei manifolds have a common Morse function
with different gluing (attaching) maps, the Morse theory might be useful for
a further study of the manifolds.

\section{Toda lattice with non-negative signs}
\renewcommand{\theequation}{2.\arabic{equation}}\setcounter{equation}{0}
\renewcommand{\thefigure}{2.\arabic{figure}}\setcounter{figure}{0}

In this paper, we use the following standard Lie theoretic notation:

\begin{Notation}
\label{standard}
Lie algebras:
\end{Notation}
Let $\frak g$ denote a real
split semisimple Lie algebra of rank $l$ with Killing form $( , )$,
universal enveloping algebra $U({\frak g})$ and we let
${\frak g}^{\prime}=Hom_{\Bbb R}({\frak g}, \Bbb R)$ denote the dual of ${\frak g}$.
Fix a Cartan decomposition $\frak g=\frak k + \frak p$ with $\frak k$ the Lie algebra
of a maximal compact Lie group $K$. We also fix a
split Cartan subalgebra ${\frak h}\subset \frak p$ with root system $\Delta $, real root vectors
$e_{\alpha_i}$ associated with simple roots $\{ \alpha_i : i=1, \cdots ,l \}=\Pi$.
Denote $\{h_{\alpha_i},e_{\pm\alpha_i}\}$ the Cartan-Chevalley basis of
$\frak g$ which satisfies the relations,
\begin{equation}
\label{chevalley}
 [h_{\alpha_i} , h_{\alpha_j}] = 0, \, \,
  [h_{\alpha_i}, e_{\pm \alpha_j}] = \pm C_{j,i}e_{\pm \alpha_j} \ , \, \,
  [e_{\alpha_i} , e_{-\alpha_j}] = \delta_{i,j}h_{\alpha_j}.
\end{equation}
where the $l\times l$ matrix
$(C_{i,j})$ is the Cartan matrix
corresponding to $\frak g$, and $C_{i,j}=\alpha_i(h_{\alpha_j})=
\langle \alpha_i,h_{\alpha_j}\rangle$.
The Weyl group $W$ is thus generated by the simple reflections $s_{\alpha_i}$, $i=1, \cdots ,l$. For any  $S\subset \Pi$, we define the subgroup generated by
$S$,
\begin{equation}
W_S=\langle ~s_{\alpha_i}~:~\alpha_i\in S~\rangle
\label{parabolic}
\end{equation}
This is the Weyl group of a parabolic Lie subgroup and it is standard to refer to $W_S$ as a {\it parabolic subgroup} of $W$.

\begin{Notation}
Lie groups:
\label{standardG}
\end{Notation}
We let $G_{\Bbb C}$ denote the connected adjoint Lie group with Lie algebra $\frak g_{\Bbb C}$
and $G$ the connected Lie subgroup correspondintg to $\frak g$.  We fix a maximal compact Lie group $U$ of $G_{\Bbb C}$ and assume that $U\cap G=K$ with Lie algebra $\mathfrak k$. 

Denote by ${\tilde G}$ the
Lie group $\left\{g\in G_{\Bbb C}~:~ Ad(g)\frak g \subset \frak g   \right\}.$
A split Cartan of ${\tilde G}$
with Lie algebra $\frak h$ will be denoted by $H_{\Bbb R}$; this Cartan subgroup has
exactly $2^l$ connected components and the component of the identity is denoted by $H$. We  let $\chi_{\alpha_i} $ denote the roots characters , defined on $H_{\Bbb R}$

\begin{Example}
\label{tildeSL}
\end{Example}

If $G=Ad(SL(n,\Bbb R) )$, ${\tilde G}$ is isomorphic to $SL(n,\Bbb R)$ for
$n$ odd and to $Ad(SL(n,\Bbb R)^{\pm}) $. This example is
the underlying Lie group for the indefinite Toda lattices as shown in 
\cite{casian:99}.  We have in this case $U=Ad(SU(n))$, $K=Ad(SO(n))$.

\subsection {Variety of ad-diagonalizable Jacobi elements and its closure \label{moment}.}



As in the introduction, we let $Z^+$ denote the set of Jacobi elements with 
positive signs in the coefficients of ${\frak g}_{\alpha}\in {\frak g}$
with $\alpha\in \pm\Pi$,
and ${\overline Z}^+$ denote its closure,
\begin{equation}
Z^{+}=\left\{X=x+\sum_{i=1}^l b_i( e_{\alpha_i}+ e_{-\alpha_i})\in{\frak g}~
:~ x\in {\frak h},~ b_i > 0\right\}
\label{bigZ}
\end{equation}

\begin{equation}
{\overline Z}^{+}=\left\{X=x+\sum_{i=1}^l b_i( e_{\alpha_i}+ e_{-\alpha_i})\in{\frak g}~
:~ x\in {\frak h},~ b_i \geq 0\right\}
\label{closedbigZ}
\end{equation}

We then define an isospectral leaf in ${\overline Z}^+$ as the
fundamental object of our study:
\begin{Definition}
 Fundamental invariants and isospectral leaves: 
\end{Definition}

Let $S(\frak g)$ be the 
symmetric algebra of $\frak g$. We may regard $S(\frak g)$ as the algebra of polynomial 
functions on the dual $\frak g^\prime$.
If we consider the algebra of $G$-invariants of $S({\frak g})$, then by Chevalley's theorem
there are homogeneous polynomials  $I_1, \cdots ,I_l$ in $S({\frak g})^{G}$ which are
algebraically independent and which generate $S({\frak g})^{G}$. According to Chevalley's
theorem, we can then express
$S({\frak g})^{G}$ as ${\Bbb R} [I_1, \cdots ,I_l]$. 
  The functions $I_1, \cdots ,I_l$ which are now polynomials on $\frak g$ can be further restricted to 
 $Z^{+}$.  This gives rise to a map function ${\cal I}=(I_1, \cdots ,I_l)$
and its restriction  ${\cal I}:Z^{+} \to {\Bbb R}^l$. We consider 
the isospectral leaf $Z^{+}(\gamma)={\cal I}^{-1}(\gamma)$
for any $\gamma \in {\Bbb R}^{I}$ in the image of ${\cal L}$. Denote by 
${\overline Z}^{+}(\gamma)$ the closure
of $Z^{+}(\gamma)$ inside ${\overline Z}^{+}$.

\vskip 0.5cm

\subsection{Moment map and convexity}
We recall the Convexity Theorem 
proven in \cite{BFR:90} for the isospectral manifold
${\overline Z}^+(\gamma)$. 
Let $\iota$ be the imbedding defined in \cite{BFR:90} 
(main theorem in section 3) of  ${\overline Z}^{+}(\gamma)$ into 
an adjoint $U$ orbit whose intersection with $\mathfrak g$ we denote  ${\cal O}_{\Lambda}$ in ${\mathfrak g}_{\mathbb C}$.

Consider the {\em moment map} associated to the $U$ orbit and restricted
to $\mathfrak g$
 $J:{\cal O}_{\Lambda}
 \to \frak h$. This is also the map obtained by orthogonal 
projection of $\frak g$ to $\frak h$ and then restriction to ${\cal O}_{\Lambda}$. Then
the composition $J\circ \iota$ has as image the convex hull of 
a Weyl group orbit
$Wx_o$ with  $x_o\in \frak h$ as in section 4.1 in \cite{BFR:90}. 
We denote the convex hull as
\begin{equation}
\label{Gamma}
\Gamma=J\circ \iota \left({\overline Z}^+(\gamma)\right)
\end{equation}
 The two objects ${\overline Z}^{+}(\gamma)$ and $\Gamma$ are {\it manifolds with boundary}. The map $J\circ \iota$
is {\it homeomorphism} according to the Convexity Theorem in \cite{BFR:90} and a diffeomorphism in the interior.

Such a convex hull is necessarily stable under the ordinary
$W$ action on $\frak h$.  Moreover there is a diffeomorphism between $H$ and the interior of  $\Gamma$ of manifolds with a $W$ action.

\subsection{Internal and external walls of $\Gamma$.}

We here devide the polytope $\Gamma$ into the Weyl chambers, and introduce the
notations to describe them.
We first denote $C^{\prime}_e$ as the 
dominant chamber in $\frak h$ intersected with $\Gamma$ and 
$\overline{C^{\prime}}_e$ the corresponding closure, and also denote 
$C_{w}^{\prime}=w(C_e^{\prime})$. We define $C_{w}=\left\{ w \right\}\times 
C^{\prime}_{w}$ and its closure $\overline {C}_{w}= \left\{ w \right\} \times \overline {C^{\prime}}_{w}$.  The $C^{\prime,\cdots}$ will refer to subsets of 
$\Gamma$, and we have the convention:
\begin{equation}
 C_{w}^{\cdots}=\left\{ w \right\}\times  C_{w}^{\prime \cdots}
\label{prime}
\end{equation}
in all our notation concerning walls.

 For each simple
 root $\alpha_i$ we may consider the corresponding 
$\alpha_i$ (internal) chamber wall intersected with $\overline {C^{\prime}}_{w}$. Denote 
this set by ${C^{\prime}}_{w}^{\alpha_i,IN}$. Each external wall of the convex 
hull of $Wx_o$ is
 parametrized by a simple roots $\alpha_i$. We denote an {\it external} 
  wall of $\Gamma$ by 
${C^{\prime}}_{w}^{\alpha_i, OUT}$ if it  intersects all the {\it internal} chamber walls except for 
${C^{\prime}}_{w}^{\alpha_i, IN}$.

For any  $A\subset \Pi$ we define the subsets of $\overline {C^{\prime}}_{w}$
of dimension $|\Pi\setminus A|$,

\begin{equation}
\begin{array} {lll}
  {C^{\prime}}_{w}^{A,\Theta} &=  \displaystyle{\bigcap_{\alpha_i\in A} } 
  {C^{\prime}}_{w}^{\alpha_i, \Theta}, 
  ~~ &{\text {if } A\not=\emptyset} ~,\\
  {C^{\prime}}_{w}^{A,\Theta} &=  C_{w}^{\prime}~, ~~  &{\text {if } A=\emptyset}~,
\end{array}
\label{walls}
\end{equation}
where $\Theta$ is either $OUT$ or $IN$. Thus we have the decomposition,
\begin{equation}
\label{cw}
{\overline {C^{\prime}}}_{w}=\displaystyle{\bigcup_{A\subset \Pi \atop
\Theta\in\{OUT,IN\}} C^{\prime A, \Theta}_{w}}.
\end{equation}

\subsection{Toda flow on $\Gamma$ and its unstable manifolds.}

Denote by $\psi_t:  \overline{Z}^+(\gamma) \to  \overline{Z}^+(\gamma)  $ the 
Toda flow on  ${Z}^+(\gamma)$ and its boundary.  For each $t$ this gives a diffeomorphism of  $\overline{Z}^+(\gamma)$ with itself  of manifolds with boundary (see \S 2 of  \cite{Milnor:65} for definitions).

The map $J\circ \iota$  between ${\overline Z}^{+}(\gamma)$ and $\Gamma$ 
allows us to obtain a Toda flow on $\Gamma$.
 Namely we have the Toda flow  $\phi_t: \Gamma \to \Gamma$ with
 \begin{equation}
 \label{todaongamma}
\phi_t=J\circ \iota \circ \psi_t \circ (J\circ \iota)^{-1}. 
\end{equation}
This flow preserves the boundary of $\Gamma$ and, the image of the critical points of $\psi_t$ are 
points in the $W$ orbit of $x_o\in C^{\prime \Pi,OUT}_e$,
the critical point in the dominant chamber.
The {\it closures} of the unstable manifolds of the Toda flow on
$\Gamma$ are easily described as follows. 
First we define the unstable Weyl group $W^u(w)$
associated to the critical point $w^{-1} x_o$ with $w\in W$,
\begin{equation}
\displaystyle{W^u(w)= W_{\Pi^u(w)}} \quad ~ {\rm with}\quad ~ \Pi^u(w):= 
\displaystyle{ \left\{\alpha_i: \ell (s_{\alpha_i} w) > \ell(w)  \right\}},
\label{tau}
\end{equation}
where $\ell(w)$ denotes the length of $w\in W$. We also define
$\Pi^s(w):=\Pi\setminus \Pi^u(w)$.

The  closure of the unstable manifold associated to the critical point 
$w^{-1} x_o$
 is then the union
\begin{equation}
{C^{\prime}}^u (w)=\displaystyle{\bigcup_{\sigma \in W^u(w) \atop
 B\supset \Pi^s(w)} }{C^{\prime}}^{B, OUT}_{(\sigma w)^{-1}}
\label{unstable}
\end{equation}

\subsection{Boundaries and cell decomposition \label{cellsTomei}}

The chamber $\overline {{{C}}^{\prime}}_{w}$ is a {\it box}. Fix two disjoint subsets $A, S \subset \Pi, A\cap S=\emptyset$. Then each intersection 
 ${C^{\prime}}^{A, OUT}_{w}\cap { C^{\prime}}^{S, IN}_{w}$ is also a box of dimension $|\Pi\setminus (A\cup S)|$. These cells can be parametrized as the triples $(A; S;  [w]_{\Pi\setminus S})$ with $[w]_{\Pi\setminus S}\in W/W_S$ the coset of $w$,
\begin{equation}
\label{intersection}
(A;S;[w] _{\Pi\setminus S})= {C^{\prime}}^{A, OUT}_{w}\bigcap { C^{\prime}}^{S, IN}_{w}.
\end{equation}
With this notation, we also have the parametrizations
\begin{equation}
\nonumber
C^{\prime A, OUT}_{w}=(A;\emptyset;w), \quad \quad C^{\prime S, IN}_{w}=(\emptyset; S; [w] _{\Pi\setminus S}).
\end{equation}

We denote the faces of $(A;S;[w]_{\Pi\setminus S})$  
by $\partial_{j,c}(A;S;[w]_{\Pi\setminus S})$ with $i=1, \cdots ,l$ 
and $c=1,2$, which are defined as
\begin{equation}
\begin{array} {ll}
& \displaystyle{  \partial_{j,1} (A; S ;[w] _{\Pi\setminus S}  )=(A\cup \left\{\alpha_j \right\};S; [w] _{\Pi\setminus S}  )} \\
&{}\\
& \displaystyle{ \partial_{j,2} (A;S ;[w] _{\Pi\setminus S})= ( A; S \cup \left\{\alpha_j \right\}; [w] _{\Pi\setminus 
( S \cup \left\{\alpha_j \right\}) })}
\end{array}
\label{boundarytriples1}
\end{equation}
\begin{figure}
\epsfysize=8cm
\centerline{\epsffile{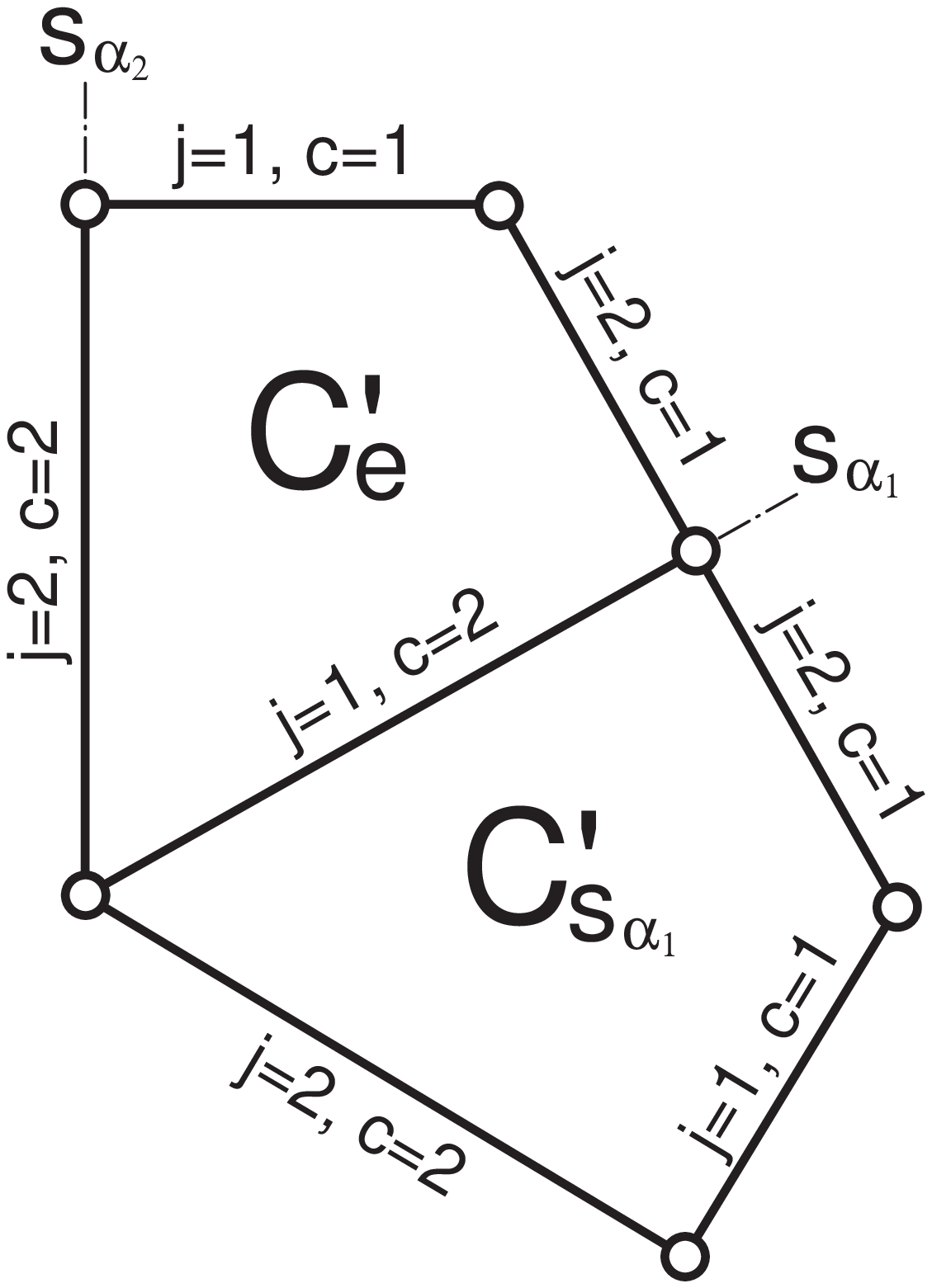}}
\caption{Boundaries $\partial_{j,c} {{{C}}^{\prime}}_{w^{-1}}$ 
for ${\frak g}$ of rank 2.}
\label{fig:c}
\end{figure}
We then define ${\cal M}_k$ the $k$-chain complex consisting of those cells,

\begin{equation}
 {\cal M}_k ={\Bbb  Z} \left[ ~(A;S; [w]_{\Pi\setminus S}) ~:~A,S \subset \Pi, ~A\cap S=\emptyset, ~
 k=|\Pi\setminus (A\cup S) |,~ w\in W ~\right]
\label{chaincomplexabstract}
\end{equation}

 To obtain cancellations along the walls of $\Gamma$  we  introduce signs  $(-1)^{\ell(w)}$ so that $\Gamma$ viewed inside ${\cal M}_l$ corresponds to the alternating sum:

\begin{equation}
\Gamma=\sum_{w\in W} (-1)^{\ell (w)} \{e\}\times{C}_{w^{-1}}
\label{gamma1}
\end{equation}

The need for the signs $(-1)^{\ell (w)}$, given our notational conventions, is illustrated in Figure \ref{fig:c}. Inner walls have the same value of $c$ on  opposite sides of a wall and then require a sign $(-1)^{\ell (w)}$. Thus cancellations along the  inner chamber walls result from the fact that two adjacent chambers correspond to Weyl group elements with lengths differing by one.  The boundary  of $\Gamma$ becomes then a cycle.

\section{ Twisted actions of $W$ }
\renewcommand{\theequation}{3.\arabic{equation}}\setcounter{equation}{0}
\renewcommand{\thefigure}{3.\arabic{figure}}\setcounter{figure}{0}

\subsection{Twisted $W$ action on ${\cal E}$}

The Weyl group associated to the Cartan subgroup $H_{\Bbb R}$ is the quotient $W=N(H_{\Bbb R})/ H_{\Bbb R}$ where $N(\cdots)$ denotes the {\it normalizer. }
Here we consider the $W$-action on the set of signs,
\begin{equation}
\label{allsigns}
\displaystyle{{\cal E}= \left\{(\epsilon_1, \cdots ,\epsilon_l)~: ~ \epsilon_i =\pm 1 \text { for all } i=1, \cdots ,l   \right\}}
\end{equation}
 with a group structure of multiplication $\bullet :{\cal E}\times{\cal E}
 \rightarrow {\cal E}$,
\begin{equation}
(\epsilon_1, \cdots ,\epsilon_l)\bullet (\epsilon^{\prime}_1, \cdots ,\epsilon^{\prime}_l)=
(\epsilon_1\epsilon_1^{\prime},\cdots ,\epsilon_l\epsilon_l^{\prime}).
\label{bullet}
\end{equation}
The group $\cal E$ parametrizes the connected components of $H_{\Bbb R}$,
and it is also in correspondence with the set of all elements $h \in H_{\Bbb R}$ such that the simple root characters  $\chi_{\alpha_i} (h)$  have absolute value one, i.e.
\begin{equation}
\label{hepsilon}
H_{\Bbb R}=\bigcup_{\epsilon\in {\cal E}} h_{\epsilon} H,  ~ \quad {\rm with} \quad 
\displaystyle{ \chi_{\alpha_i} (h_\epsilon ) = \epsilon_i}.
\end{equation}

Let us now define:
\begin{Definition}
 Twisted actions of the Weyl group on ${\cal E}$ \label{twisteddef2}:
\end{Definition} 

We call  an action  of $W$ on  ${\cal E}$ a {\it twisted sign action } if
there is a matrix of integers $a_{i,j}$ with $i,j \in \left\{1, \cdots ,l \right\}$ and $a_{i,i}=2$ (or just even) such that
for each $\alpha_i$ the action $s_{\alpha_i}\in W: {\cal E} 
\rightarrow {\cal E}$ is given by
\begin{equation}
 \displaystyle{s_{\alpha_i}:\epsilon_j \mapsto 
\epsilon_j\epsilon_i^{a_{j,i}}}.
\label{smatrix2}
\end{equation}

Two canonical
examples are the usual $W$ action on $H_{\Bbb R}/H$  (and thus  on $\cal E$) where $a_{i,j}=-C_{i,j}$ the Cartan matrix (recall
$s_{\alpha_i}\chi_{\alpha_j}=\chi_{\alpha_j}\chi_{\alpha_i}^{-C_{j,i}}$), 
and the trivial $W$ action where $a_{i,j}$ is even for all $i,j$.

\medskip

By (\ref{smatrix2})  the action of $s_{\alpha_{i}}$ on any $\epsilon=(\epsilon_1, \cdots ,\epsilon_l)$ such that $\epsilon_i=1$ is trivial. Hence the  only relevant cases happen when $\epsilon_i=-1$.  
For each pair of simple roots $\alpha_i, \alpha_j$  $i < j$  which are connected in the Dynkin diagram, we define
\begin{equation}
s_{ij}=\epsilon_i^{a_{j,i}}=(-1)^{a_{j,i}}
\label{signss}
\end{equation}
Note that   (\ref{smatrix2}) ensures that  any $W$  action  in Definition (\ref{twisteddef2}) necessarily acts on $\cal E$ by group automorphisms. To see this we write the $j$ th components of $s_{\alpha_i}(\epsilon\bullet \epsilon^{\prime})$ and  $s_{\alpha_i}\epsilon \bullet s_{\alpha_i}\epsilon^{\prime}$; these  are respectively $(\epsilon_i\epsilon^{\prime}_i)^{a_{j,i}}\epsilon_j\epsilon^{\prime}_j$
and $\epsilon_i^{a_{j,i}}\epsilon_j (\epsilon^{\prime})^{a_{j,i}}\epsilon^{\prime}_j$.

The twisted sign actions can be defined by first modifying the usual action of $W$ on a split Cartan subgroup. The following definition describes the kind of modification that leads to Definition (\ref{twisteddef2}):

\begin{Definition}
 Twisted action of the Weyl group on $H_{\Bbb R}$\label{twisteddef}:
\end{Definition}

Denote by $*$ the usual action of $W$ on the Cartan subgroup $H_{\Bbb R}$.
An  action of $W$ on $H_{\Bbb R}$ is said to be a {\it twisted action } if it satisfies

\begin{enumerate}
\item $W$ acts on $H$ with the usual action; moreover if $\chi_{\alpha_i} (h)$ is positive,
 then $s_{\alpha_i} h=s_{\alpha_i}* h$ \label{blue} 
 \item  each $w\in W$ acts as a Lie group automorphism of $H_{\Bbb R}$ and thus
$w(h_1h_2)=w(h_1)w(h_2)$ for any pair of elements $h_1, h_2$ in $H_{\Bbb R}$.\label{auto}
\end{enumerate}

\medskip

Note that any twisted action determines an action on the set of connected
components $H_{\Bbb R}/ H$ which is in bijective correspondence with the set of  signs ${\cal E}$.  The actions obtained in this way , starting from Definition (\ref{twisteddef}), give rise to the twisted sign actions.

\subsection{The rank 2 cases \label{rank2cases}}

We abuse notation slightly denoting by $s_{\alpha_i}$ what should be $A(s_{\alpha_i})$ for an appropriate group homomorphism $A:W\to Aut{\cal E}$ . In these cases, an easy calculation shows (see also below) that  the pair of signs $(s_{12}, s_{21})$ always determine constructions which satisfy relations of the following kinds

\begin{eqnarray}
\begin{array}{ccc}
& (s_{\alpha_i}s_{\alpha_j})^2=e\quad
& s_{12}\not=s_{21} \\
&{}&{} \\
& (s_{\alpha_i}s_{\alpha_j})^3=e\quad
& s_{12}=s_{21}=-1 \\
\label{puf}
\end{array}
\end{eqnarray}

By squaring the first relation one obtains $(s_{\alpha_i}s_{\alpha_j})^4=e$ and
by squaring the second one obtains $(s_{\alpha_i}s_{\alpha_j})^6=e$.  Thus the
first case is compatible with the braid relation of $W$ in $B_2$ and $G_2$ and the second is compatible with the cases of
$A_2$ and $G_2$.  We thus list the twisted sign actions in all these cases.

\vskip 0.5cm
\begin{Example}
\label{Relations}
\end{Example}
We just illustrate what is involved in checking  (\ref{puf}) by working out the case, $s_{12}=-1$ and $s_{21}=1$.
For example, if we apply $(s_{\alpha_2}s_{\alpha_1})^2$ to $(-1,-1)$ we get again $(-1,-1)$ as demonstrated below:

\begin{equation}
(-1,-1)\overset{s_{\alpha_1}}{\longrightarrow}(-1,1)\overset{s_{\alpha_2}}
{\longrightarrow}
(-1,1)\overset{s_{\alpha_1}}{\longrightarrow}(-1,-1)\overset{s_{\alpha_2}}
{\longrightarrow}(-1,-1).
\label{B2rel}
\end{equation}
Similarly we obtain  that $(s_{\alpha_2}s_{\alpha_1})^2$ gives the identity by considering all the other cases. Hence the relation $(s_{\alpha_2}s_{\alpha_1})^2=e$ is obtained for
$s_{12}\ne s_{21}$ as in (\ref{puf}).

\begin{Example}
The twisted sign actions in $A_2$:
\label{A2case}
\end{Example}

The only non-trivial twisted sign action corresponds to $s_{12}=-1$ and $s_{21}=-1$, 
the standard sign action, and it is given explicitly by

\begin{equation}
\begin{array} {ll}
& \displaystyle{s_{\alpha_1}(\epsilon_1, \epsilon_2)=(\epsilon_1, \epsilon_1\epsilon_2)} \\
& \displaystyle{s_{\alpha_2}(\epsilon_1, \epsilon_2)=(\epsilon_2\epsilon_1,\epsilon_2)} \\
\end{array}
\label{A2action}
\end{equation}

\vskip 0.5cm

\begin{Example}
The twisted sign actions of $B_2$ and $C_2$:
\label{B2case}
\end{Example}

The only non-trivial possibilities here are
the cases $s_{12}\not=s_{21}$. These become explicitly:

\begin{equation}
\begin{array} {ll}
& \displaystyle{s_{\alpha_1}(\epsilon_1, \epsilon_2)=(\epsilon_1, \epsilon_1\epsilon_2)}, \\
& \displaystyle{s_{\alpha_2}(\epsilon_1, \epsilon_2)=(\epsilon_1,\epsilon_2)} ,
\end{array}
\label{B2action}
\end{equation}
which is the standard action for $B_2$ and a twisted one for $C_2$. We also have
\begin{equation}
\begin{array}{ll}
& \displaystyle{s_{\alpha_1}(\epsilon_1, \epsilon_2)=(\epsilon_1, \epsilon_2)}, \\
& \displaystyle{s_{\alpha_2}(\epsilon_1, \epsilon_2)=(\epsilon_2\epsilon_1,\epsilon_2)}. \\
\end{array}
\label{C2action}
\end{equation}
for the standard case of $C_2$ and a twisted one of $B_2$.
\vskip 0.5cm

\begin{Example}
The twisted sign actions in $G_2$:
\label{G2case}
\end{Example}

Assume
that $C_{21}=3$. In this case there are four possible
specifications of the signs $s_{12}$ and $s_{21}$ and they all give
twisted sign actions. The standard sign action corresponds to $s_{12}=-1$ and $s_{21}=-1$.
\vskip 0.5cm

\begin{Definition}
Marked Dynkin diagrams
\label{marked}
\end{Definition}
 A {\it marked} Dynkin diagram is a Dynkin diagram where all the single edges are assigned
 the sign $+$ or $-$; any double edge is assigned a pair of signs $+,+$
or $\pm,\mp$, and a triple edge ($G_2$) is assigned any pair of any signs $\pm,\pm$ or $\pm,\mp$.  Those signs associated to the edge joining
$\alpha_i$ and $\alpha_j$ are given by $s_{ij}$ or $s_{ji}$, and  $s_{ij}=s_{ji}$ in the case of single edges in the Dynkin diagram.

\medskip

We call a marked Dynkin diagram {\it positively marked} if all the signs $s_{ij}$ over any of its edges are positive (the case  when $\frak g$ is of type $A_1\times \cdots \times A_1$ is always considered positive). A positively marked Dynkin diagram corresponds to  a trivial action of $W$ on $\cal E$. We call the marked Dynkin diagram {\it standard} in case that the matrix $a_{ij}$ can be chosen to be the Cartan matrix.  Given a marked Dynkin diagram and 
subset $B\subset \Pi$, there is a subdiagram $\delta_B$ corresponding to an action of $W_B$ on a set of signs. 

\begin{Proposition} Assume that $\frak g$ is a simple Lie algebra. Let 
$\alpha$
denote the number of pairs of simple roots joined by a single edge in the Dynkin
diagram and $\beta$ the number of double edges. 
Then if $\frak g$ 
is not of type $G_2$  there are exactly $2^{\alpha}3^{\beta}$ twisted 
signed actions of $W$ parametrized by marked Dynkin diagrams. 
In the case of $G_2$ there are four such twisted signed actions. For any 
semisimple Lie algebra the twisted signed actions are parametrized by the
set of marked Dynkin diagrams. \label{classificationtwisted}
  
\end{Proposition} 

\begin{Proof} Recall that $W$ is a Coxeter group, (Proposition 3.13  \cite{kac:90}) and it
thus has defining relations $s_{\alpha_i}^2=e$ and $(s_{\alpha_i}s_{\alpha_j})^{m_{ij}}=e$ where 
 $m_{ij}$ is $2,3,4,6$ depending on the number of lines joining $\alpha_i$ and $\alpha_j$
 in the Dynkin diagram.  The case $m_{ij}=2$ occurring when $\alpha_{i}$ and $\alpha_{j}$  are not connected in the Dynkin diagram.  From the definition of a twisted signed action of $W$ on $\cal E$,
 such an action  is completely determined if one
 specifies all the signs $s_{ji}:=(s_{\alpha_i}\epsilon )_j$ 
 where $\alpha_j$ is joined
 in the Dynkin diagram to $\alpha_i$. It is then enough to count all the possible choices
 $s_{ij}$ and $s_{ji}$ giving rise to 
 $(s_{\alpha_i}s_{\alpha_j})^{m_{ij}}=e$; the other relation being automatically satisfied.
 This then reduces the proof to the classification of all twisted signed $W$ actions in the cases
 $A_2$, $B_2$, $G_2$, namely the  examples (\ref{A2case}), (\ref{B2case}) and (\ref{G2case}) in subsection (\ref{rank2cases}).
 To a marked Dynkin diagram we thus associate a twisted sign action of $W$ where the signs
 $s_{i,j}$ defining the action are given as in definition (\ref{marked}).
\end{Proof}

\section{The twisted Tomei manifolds}
\renewcommand{\theequation}{4.\arabic{equation}}\setcounter{equation}{0}
\renewcommand{\thefigure}{4.\arabic{figure}}\setcounter{figure}{0}

For each marked Dynkin diagram $\delta$ corresponding to a twisted sign  action we associate a compact smooth manifold with an action of $W$. 
We call the manifold the twisted Tomei manifold, and here give
a consruction of the manifold by gluing the chamber walls of the 
polytopes associated with different $ \epsilon \in {\cal E}$:

 Let $\delta$ denote a marked Dynkin diagram giving rise to a
 fixed twisted sign action of $W$. 
 Recall {\it external } and {\it internal} chamber walls
which constitute the boundary of $\Gamma$ and were described above.

We then define gluing maps between the chamber walls
denoted by $\{\epsilon\}\times C_{w}^{\cdots}=\{\epsilon\}\times
\{w\}\times C_{w}^{\prime \cdots}$ as follows:
For the internal walls, we define
\begin{equation}
\begin{array} {lll}
g_{w,i,IN} : &\displaystyle{\left\{\epsilon \right\}\times 
 C_{w}^{\alpha_i, IN} \to \left\{s_{\alpha_i}\epsilon \right\}\times 
 C_{ w s_{\alpha_i}}^{\alpha_i, IN}     } \\
& ~~\displaystyle{(\epsilon ,w, x ) \longmapsto (s_{\alpha_i}\epsilon , ws_{\alpha_i} , x) }
\end{array}
\label{glueIN}
\end{equation}
where note $ws_{\alpha_i}w^{-1}x=x$.
For the external walls, we define
\begin{equation}
\begin{array} {ll}
  g_{w,i,OUT}: &\displaystyle{\left\{\epsilon \right\}\times C_{w}^{\alpha_i, OUT} \to \left\{\epsilon^{(i)} \right\}\times 
C_{w}^{\alpha_i,OUT}     } \\
& ~~\displaystyle{(\epsilon , w, x ) \longmapsto (\epsilon^{(i)} , w , x) }
\end{array}
\label{glueOUT}
\end{equation}
where $ \epsilon^{(i)}=(\epsilon_1, \cdots ,-\epsilon_i, \cdots ,\epsilon_l)$.
We denote ${\tilde M}(\delta)$ the disjoint union of all the chambers with 
different signs,
\begin{equation}
\tilde {M }(\delta)=\bigcup_{w\in W \atop
\epsilon \in {\cal E}} \{w\epsilon\}\times \overline {C}_{w^{-1}} .
\label{disjoint}
\end{equation}
We also denote $M(\delta)$ the topological space 
obtained from the disjoint union
in   (\ref{disjoint}) by gluing along the internal and
external walls  using the maps
$g_{w,i,IN}$ and $g_{w,i,OUT}$. There is then a map 
\begin{equation}
z:\tilde{M}(\delta)\to M(\delta).
\label{z}
\end{equation}

\begin{figure}
\label{fig:gluing}
\epsfysize=7.5cm
\centerline{\epsffile{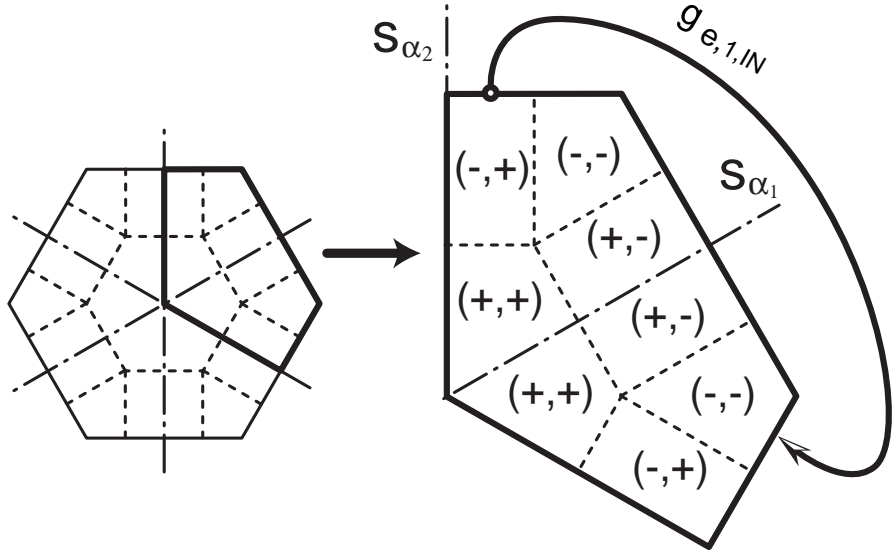}}
\caption{Standard case $A_2$; gluing by
$g_{e,1,IN}$ along
the internal walls of $\{(-,+)\}\times
C^{\alpha_1,IN}_e$ and $\{(-,-)\}\times C^{\alpha_1,IN}_{s_{\alpha_1}}$
reverses orientation (internal hexagon is $\Gamma$).
External walls are glued along the dashed lines. }
\end{figure}
We define an action of the group $W$ on ${\tilde M}(\delta)$ by
\begin{equation}
\sigma\in W : (\epsilon , w, x ) \mapsto
(\epsilon , \sigma w , \sigma x ).
\label{extendedaction}
\end{equation}

\begin{Lemma}The action of $W$  defined in (\ref{extendedaction}) gives rise to a well defined action in $M(\delta)$\label{actionlemma}.

\end{Lemma}
\begin{Proof}
We must show that the action of $W$ sends pairs of  elements that are identified with each other to other pairs of elements that are identified.
Consider first the case of the gluing maps for internal walls $g_{\cdots,i,IN}$:
Let the vertical arrows in  (\ref{checkaction}) below correspond to multiplication by $\sigma\in W$;   and assume the horizontal arrows are given by a gluing map of the form $g_{w, i,IN}.$   We have:

\begin{equation}
\begin{CD}
(\epsilon, w, x) @>{g_{w, i,IN}}>> (s_{\alpha_i}\epsilon , ws_{\alpha_i}, x)\\
@V{\sigma}VV  @V{\sigma}VV\\ 
(\epsilon , \sigma w, \sigma x) @>{g_{\sigma w,i, IN}}>> (s_{\alpha_i} \epsilon , \sigma ws_{\alpha_i}, \sigma x)
\end{CD}
\label{checkaction}
\end{equation}
We thus
 obtain that the image of the first pair of points under the group action 
corresponds to another pair of points which are identified. The gluings along the
external walls do not pose any difficulties either.
\end{Proof}

\begin{Remark}
\end{Remark}
Consider the union of the chambers with the signs determined by a twisted action of W corresponding to $\delta$,
 
 \begin{equation}
\displaystyle{ \tilde{\Gamma}_{\epsilon}(\delta)=\bigcup_{w\in W} \left\{ w\epsilon \right\}\times \overline {C}_{w^{-1}}}
\label{gamaw}
\end{equation}
We let  $\Gamma_{\epsilon}(\delta)$ denote the image of $\tilde{\Gamma}_{\epsilon}$ in $M(\delta)$. Then
 after the identifications in $M(\delta)$, this space becomes a copy of $\Gamma$ (see Figure \ref{fig:0}). We denote 

\begin{equation}
F_{\epsilon}^{\delta} ~:~\Gamma\longrightarrow \Gamma_\epsilon(\delta)
\label{Fepsilon}
\end{equation}
We note that each $\Gamma_{\epsilon}$ is usually {\it not} preserved by the $W$ action. The union of  the interiors of all the $\Gamma_{\epsilon}$ over
$\epsilon \in {\cal E}$ can be identified with a copy of $H_{\Bbb R}$  and the $W$ action would then correspond to the notion of a twisted $W$ action on a Cartan subgroup. We can now express $M(\delta)$ as a union:

\begin{equation}
M(\delta)=\bigcup_{\epsilon\in {\cal E}} \Gamma_\epsilon(\delta).
\end{equation}
 
\vskip 0.5cm
\noindent
 Similarly we consider
 \begin{equation}
\label{chamber}
\displaystyle{ \tilde {K}_e=\bigcup_{\epsilon\in \cal {E}} \left\{ \epsilon \right\}\times \overline{C}_e   } 
\end{equation}
 Denote by $K_e=z(\tilde {K}_e)$ its image in $M(\delta)$.  This set is compact and constitutes a fundamental domain of the $W$ action. Moreover $K_e$ is a {\it box}. For example  Figure \ref{fig:1} shows the union of $K_e$ and $s_{\alpha_1}K_e$ in the $A_2$ case.

\begin{Proposition} The space $M(\delta)$ has the structure of a smooth compact manifold with
a $W$ action.\label{ourmanifold}
\end{Proposition}

\begin{Proof} The space $M(\delta)$ has an action of  $W$ as defined in  (\ref{extendedaction}) and verified in  (\ref{checkaction}) in Lemma \ref{actionlemma}. 
The compactness follows from the fact that  $W$
is finite and that the box $K_e$ is compact because $W(K_e)=M(\delta)$. The smoothness is obtained using the coordinate charts given by the boxes $\{w(K_e): w\in W\}$. 
\end{Proof}

\subsection{Toda flow on $M(\delta)$}

Using the map $F_\epsilon^{\delta}$ in (\ref{Fepsilon}),
we define the Toda flow

\begin{equation}
\begin{array} {ll}
\phi^{\delta}_t : & \displaystyle{ M(\delta)\longrightarrow M(\delta)} \\
& \displaystyle{(\epsilon, w, x)\mapsto F_{\epsilon}^{\delta}\circ\phi_t\circ (F^{\delta}_{\epsilon})^{-1}
 (\epsilon, w,x)}
\end{array}
\label{Todaflow}
\end{equation}
where $\phi_t$ is the Toda flow on $\Gamma$ defind in (\ref{todaongamma}).

Since the unstable manifolds for the Toda flow on $\Gamma$ have been described in  (\ref{unstable}) we are now able to describe the {\it closures}  $\overline{C}^u(w)$ of the unstable manifolds $C^u(w)$ of the new Toda flow in the twisted Tomei manifold.

Recall the map $z : {\tilde M}(\delta)\to M(\delta)$ in  (\ref{z}). Then 
we have
\begin{equation}
\overline {C}^u(w)=\bigcup_{\sigma \in W^u(w) \atop B\supset \Pi^s(w)}z
( {\cal E} \times  {C}^{B, OUT}_{(\sigma w)^{-1}})
\label{unstable2}
\end{equation}
Since the closures of unstable manifolds are constructed in the same way as the twisted Tomei manifolds with respect to a convex polytope associated to a parabolic subgroup of $W$, we have:

\begin{Proposition} For any $w\in W$, the closure of the unstable manifold 
$\overline {C}^u(w)$ is smooth.  Moreover $\overline {C}^u(w)$ is a twisted
Tomei manifold for  the Levi factor of a parabolic subgroup determined by  $\Pi^u(w)$ with marked Dynkin diagram $\delta_{\Pi^u(w)}$ (see Definition \ref{marked}).\label{smoothunstable}
\end{Proposition}

\subsection{Cell decomposition  of $M(\delta)$ }

We now describe a cell decomposition that is useful in our explicit description of  all the unstable manifolds of the Toda lattice. This cell decomposition is not  efficient in terms of computing homology because it contains too many cells; however  a more efficient cell decomposition is described in  \cite{casian:99} in the standard case. Unfortunately  the Toda lattice is not  a Morse-Smale vector field and consequently the boundary of an unstable manifold  usually fails to be a combination of unstable manifolds. This problem  is evident in  Lemma \ref{twos} or in low dimensional examples (see also
\cite{casian:00}).

Using the cell decomposition in subsection \ref {cellsTomei} we obtain a cell decomposition for $M(\delta)$ by simply adding a factor $\left\{ \epsilon \right\} \times (\cdots)$  and applying $z$.
Note here that because of the gluings between external walls  in (\ref{glueOUT}), the cells $z(\left\{ \epsilon \right\}
 \times {C^{\prime}}^{A, OUT}_{w}\cap { C^{\prime}}^{S, IN}_{w})$  are not parametrized simply by the triples involving $\epsilon\in {\cal E}$ and $A,S\subset \Pi$. The value of $\epsilon_i$ for $\alpha_i\in A$ will be set to zero to avoid duplication of parameters. For any $A\subset \Pi$, we introduce

\begin{equation}
{\cal E}^A=\Big\{(\epsilon_1,\cdots,\epsilon_l): \epsilon_i=\pm 1 \text { if  }\alpha_i\not\in A, \epsilon_i=0  \text { if }  \alpha_i\in A\Big\}
\label{newE}
\end{equation}
Note that $W_{\Pi\setminus A}$ acts on ${\cal E}^A$, the new zeros are irrelevant to the action.

\begin{Definition}
Signed colored Dynkin diagrams 
\label{signedcol}
\end{Definition}
A {\it signed colored} Dynkin diagram is a triple $(\epsilon; A;S)$, where $\epsilon=(\epsilon_1,\cdots\epsilon_l)$ with $\epsilon_i=0$ if $\alpha_i\in A$ (otherwise $\epsilon_i=\pm1$)  and $A,S\subset \Pi$ with  $A\cap S=\emptyset$. This corresponds to the notion of  signed colored Dynkin diagram in Definition (5.1.2) in  \cite{casian:99}. There the simple roots in $S$ are colored  {\it blue } if $\epsilon_i=1$, {\it red } if $\epsilon_i=-1$; and are given a sign  $=sign \epsilon_i=\pm$ if the simple root is in $\Pi\setminus (A\cup S)$ and a zero if the simple root belongs to $A$.

\begin{figure}
\epsfysize=8cm
\centerline{\epsffile{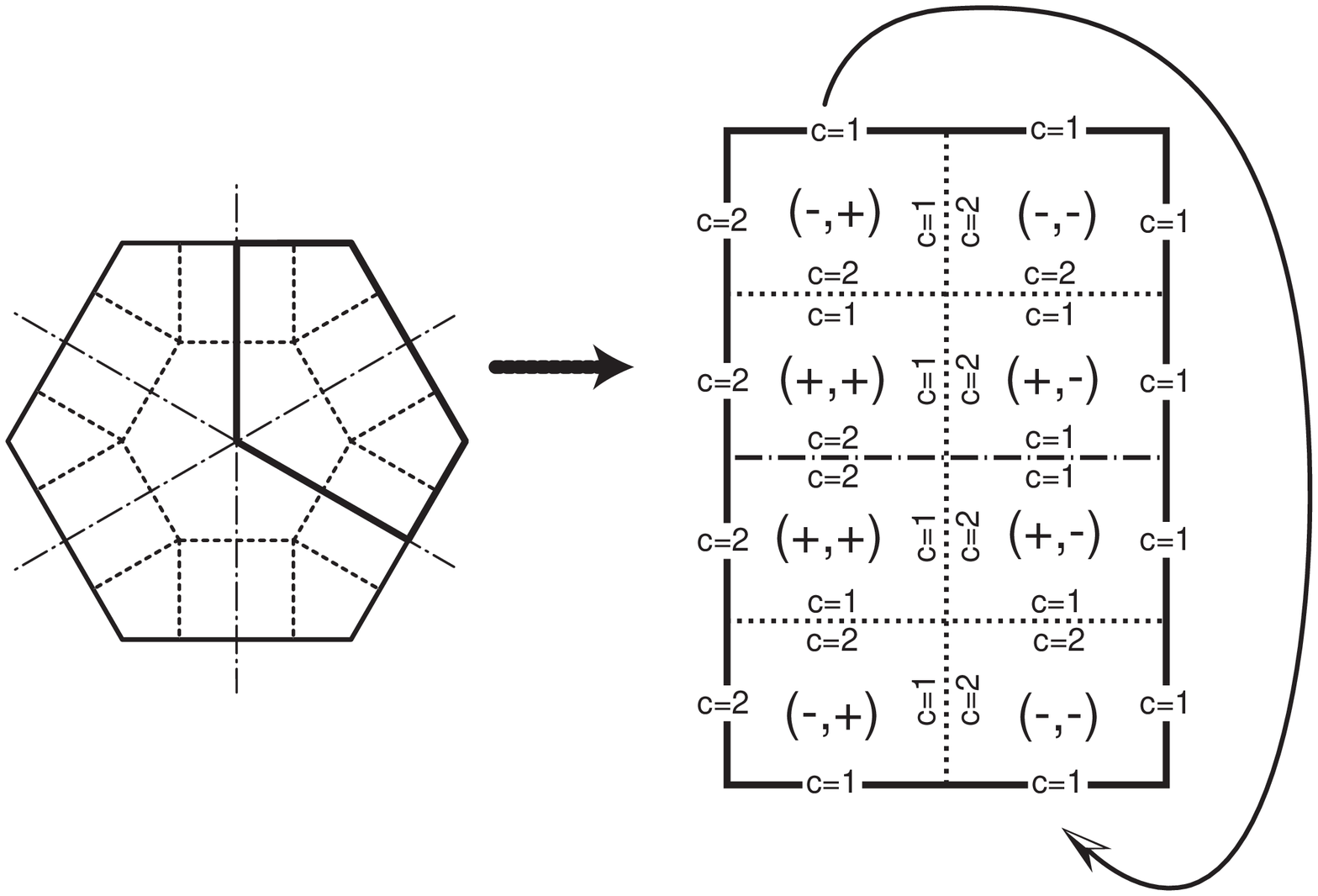}}
\caption{Boundaries $\partial_{j,c}$ for $A_2$. Cancellations along
the external walls are shown by the different choice of $c$.}
\label{fig:1}
\end{figure}

Then the cell $ \left\{ \epsilon \right\}
 \times {C^{\prime}}^{A, OUT}_{w}\cap { C^{\prime}}^{S, IN}_{w}$ is parametrized by a signed colored Dynkin diagram,
\begin{equation}
\label{scdynkinforintersection}
\left\{ \epsilon \right\}
 \times {C^{\prime}}^{A, OUT}_{w}\cap { C^{\prime}}^{S, IN}_{w}
 =(\epsilon;A;S;[w]_{\Pi\setminus S})
 \end{equation}
The faces of this cell are then expressed as 
\begin{equation}
\begin{array} {lll}
& \displaystyle{  \partial_{j,c} (\epsilon;A;S;[w]_{\Pi\setminus S}  )
=(\epsilon;A\cup \left\{\alpha_j \right\};S; [w]_{\Pi\setminus S}  )}, \quad   &\text { if }  (-1)^{c+1}=\epsilon_j\\
&{}\\
& \displaystyle{ \partial_{j,c} (\epsilon;A;S;[w]_{\Pi\setminus S})= (
\epsilon;A; S \cup \left\{\alpha_j \right\}; [w]_{\Pi\setminus ( S \cup \left\{\alpha_j \right\}) }) },\quad &\text { if }  (-1)^{c}=\epsilon_j\\
\end{array}
\label{boundarytriples2}
\end{equation}
This is  just  (\ref{boundarytriples1}) but  all the $2^l$ signs are now used to form a bigger box as in Figure (\ref{fig:1}) and cancellations are required along the {\it external} walls (which have become internal in the box).

We now define the chain complex ${\cal M}_*$,
\begin{equation}
\label{chaincomplexM}
0\longrightarrow {\cal M}_l\overset {\partial_l}{\longrightarrow}
{\cal M}_{l-1}\overset{\partial_{l-1}}{\longrightarrow}
\cdots\overset{\partial_2}{\longrightarrow}{\cal M}_1
\overset{\partial_1}{\longrightarrow}{\cal M}_0\longrightarrow 0,
\end{equation}
where the $k$-chain ${\cal M}_k$ and the boundary map
$\partial_k$ are defined as
\begin{equation}
 {\cal M}_k ={\Bbb  Z} \left[ (\epsilon;A;S; [w]_{\Pi\setminus S})~:~
  A,S\subset \Pi, A\cap S=\emptyset,  k=|\Pi\setminus (A\cup S) |, w\in W,  \epsilon\in {\cal E} \right],
\label{chaincomplex2}
\end{equation}
\begin{equation}
\partial_k=\sum_{j=1,\cdots ,k \atop  c=1,2}(-1)^{j+c+1} \partial_{j,c}.
\label{boundary3}
\end{equation}

\subsection{The boundaries of the unstable manifolds in $M(\delta)$.}
With the gluing along the walls, one needs to get the cancellations
of the orientations of the glued chambers. This then leads to a change
of orientation at the boundary.
 In order to obtain cancellations along the internal walls of $\Gamma_{\epsilon}$  we  had introduced signs  $(-1)^{\ell(w)}$ in (\ref{gamma1}). Recall (also see Figure \ref{fig:1}) that external wall cancellations occur because $c$ differs on two sides of the wall,  and internal wall cancellations  are produced by the different signs $(-1)^{\ell (w)}$.  The exceptions to this are all the cases of non-trivial gluings that change the orientation; these orientation changes can only happen along the internal  $\alpha_i$ walls satisfying $\epsilon_i=-1$ ({\it  negative walls}) . In these cases additional signs must be introduced to compensate orientation changes.  For fixed $w$ and $\epsilon$ we denote $(w\epsilon )_j$ the $j$-th component of $w\epsilon$. The cell $\Gamma_\epsilon$ can then be represented in ${\mathcal M}_*$ by

\begin{equation}
\Gamma_{\epsilon}=\sum_{w\in W} (-1)^{\ell (w)} \left( { \prod_{k=1,\cdots ,l}\epsilon_k  (w{\epsilon} )_{k}} \right)  \left\{w \epsilon \right\}\times {C}_{w^{-1}}
\label{gamma2}
\end{equation}
The boundary  of each $\Gamma_\epsilon$ becomes then a cycle and only contains external walls.

\begin{Example} The case of $A_n$.

\end{Example}
We illustrate the derivation of the expression in (\ref{gamma2}) in the case of $A_n$.
Canonical examples are the original Tomei manifold and the standard one 
($\delta$ is positively marked  or  $\delta$ is marked with signs determined by the Cartan matrix as in Definition \ref{marked}). The change of orientation of the boundary of $K_e$ under the gluing  $g_{w,i,IN}$ in (\ref{glueIN}), with $i=1$ is determined by the number of $\epsilon_j$ that change sign under the action of $s_{\alpha_1}$, that is by  the sign of $\prod_{k=1,\cdots l} \epsilon_k (s_{\alpha_1}\epsilon )_k$. This number is $1$ in the positively marked case and $-1$ in the standard case.  Therefore, to have cancellations at the internal wall boundaries we have to insert a sign $\prod_{k=1,\cdots l} \epsilon_k (s_{\alpha_1}\epsilon )_k$.  If we now consider the $\alpha_2$ wall of the new chamber $\left\{s_{\alpha_1}\epsilon \right\}\times C_{s_{\alpha_1}}$, the new sign $s_{\alpha_1}\epsilon$  gets changed under the gluing map to $s_{\alpha_2}s_{\alpha_1}\epsilon$. Thus the change of sign is given by  the product $\prod_{k=1, \cdots, l}  (s_{\alpha_2}s_{\alpha_1}\epsilon)_k (s_{\alpha_1}\epsilon )_k$. However there is already a sign attached to this chamber, namely  $\prod_{k=1,\cdots,l} \epsilon_k (s_{\alpha_1}\epsilon )_k$ and the product of these two gives the new sign which becomes $\prod_{k=1, \cdots, l} \epsilon_k (w\epsilon)_k$ with $w=s_{\alpha_2}s_{\alpha_1}$. This leads to  the expression in (\ref{gamma2}).

\vskip 0.5cm
The unstable manifold associated with $w\in W$
is now represented in ${\cal M}_*$ by

\begin{equation}
C^u(w)=\displaystyle{\sum_{\sigma \in W^u(w) \atop \epsilon \in {\cal E}^{\Pi^s(w)}} (-1)^{\ell (\sigma)+\ell (w)}} \left( \prod_{\alpha_k\in \Pi^u(w)} \epsilon_k (\sigma\epsilon)_k \right) \left\{ \epsilon \right\} \times (\epsilon;\Pi^s(w);\emptyset; (\sigma w)^{-1}).
\label{unstableclassdynkin}
\end{equation}
Here recall that $(\epsilon;\Pi^s(w);\emptyset; 
(\sigma w)^{-1})=C^{\Pi^s(w),OUT}_{(\sigma w)^{-1}}$.
The dimension of $C^u(w)$ is then given by $|\Pi^u(w)|$, and we call
it the {\it index} of $w\in W$.
Then we have:

\begin{Lemma} Let $w\in W$ of index $k$. Then there is $X\in {\cal M}_{k-1}$ such that $\partial C^u(w)=2X$. Explicitly we have
\begin{equation}
\partial C^u(w)=2(-1)^{\ell (w)+1}\displaystyle{\sum_{\sigma\in W^u(w),~
\alpha_r\in \Pi^u(w) \atop \epsilon\in {\cal E}^{\Pi^s(w)\cup \left\{\alpha_{r}\right\}  }}}
\mu(\epsilon, r, \sigma) (\epsilon; \Pi^s(w)\cup \left\{\alpha_{r}\right\};
\emptyset; (\sigma w)^{-1})
\label{explicit}
\end{equation}
where the coefficient $ \mu(\epsilon, r, \sigma)$ is given as
\begin{equation}
\left\{
\begin{array} {ll}
{ (-1)^{\ell (\sigma) + r}\displaystyle{\prod_{\alpha_j\in \Pi^u(w)} }\epsilon (-r)_j (\sigma\epsilon (-r))_j},& { \text{ if } \displaystyle{\prod_{\alpha_j\in \Pi^u(w)}}(\sigma\epsilon (-r))_j(\sigma\epsilon(+r))_j = 1} \\
&{}\\
{ ~~~~~~~~~~~0},& {\text{ if } \displaystyle{\prod_{\alpha_j\in \Pi^u(w)} }(\sigma\epsilon (-r))_j (\sigma\epsilon (+r))_j=-1}
\end{array}
\right.
\label{two3}
\end{equation}
where $\epsilon(\pm r)=(\epsilon_1,\cdots,\overset{r}{\pm 1},\cdots,\epsilon_l)$
for $\epsilon\in {\cal E}^{\Pi^s(w)\cup\{\alpha_r\}}$.
 \label{twos}
\end{Lemma}
\begin{Proof} We apply the definition of the boundary $\partial$ in (\ref{boundarytriples2}) and  (\ref {boundary3})   to the expression  in (\ref{unstableclassdynkin}). Then  noticing $\epsilon_r=(-1)^{c+1}$ we obtain

\begin{equation}
\partial C^u(w)=  (-1)^{\ell (w)} \displaystyle{
\sum_{\sigma \in W^u(w), ~\alpha_r \in \Pi^u(w)\atop
\epsilon \in {\cal E}^{\Pi^s(w)} }}  
 \nu(\epsilon,r,\sigma)
  (\epsilon[r];\Pi^s(w)\cup\left\{\alpha_{r}\right\} ; \emptyset ; 
  (\sigma w)^{-1})
\label{unstableclassdynkin2}
\end{equation}
where $\epsilon[j]=(\epsilon_1,\cdots,\overset{r}{0},\cdots,
\epsilon_l)$ for $\epsilon\in {\cal E}^{\Pi^s(w)}$, and
$\nu(\epsilon,r,\sigma)$ is given by
\begin{equation}
\label{nu}
\nu(\epsilon,r,\sigma)= (-1)^{\ell (\sigma)+r}\epsilon_r \displaystyle{\prod_{\alpha_j\in \Pi^u(w)} \epsilon_j (\sigma\epsilon)_j } 
 \end{equation}

We  will now simply collect terms with the same $\alpha_r$ but different sign $\epsilon_r$.  We then fix $\epsilon \in {\cal E}^{\Pi^s(w)\cup\left\{\alpha_{r}\right\}}$ so that 
$\epsilon (\pm r)\in {\cal E}^{\Pi^u(w)}$ with different signs in the $r$-th component.  The coefficient of   $(\epsilon;  \Pi^s(w)\cup \left\{\alpha_{r}\right\};\emptyset; (\sigma w)^{-1})$ becomes (with $\epsilon (-r)_r=-1$ by convention)

\begin{equation}
(-1)^{\ell (w)+r}(\epsilon (-r)_r)( \prod_{ j\not=r} \epsilon (-r)_j)\left( \epsilon (-r)_r\prod_{k} (\sigma\epsilon (-r))_k  +\epsilon (-r)_r \epsilon (r)_r^2 \prod_{k} (\sigma\epsilon (r)_k)\right)
\label{two1}
\end{equation}
From (\ref{two1}) and since $\epsilon (-r)_r=-1$ and $\epsilon (r)_r=1$ , the only cases when a coefficient is zero are those for which:  
\begin{equation}
  \prod_{k} (\sigma\epsilon (-r))_k= -\prod_{k} (\sigma\epsilon (r))_k
\label{two2}
\end{equation}
We thus obtain  (\ref{two3}) for the coefficients $ \mu(\epsilon, r, \sigma)$.
From here it follows that  $\partial C^u(w)=2X$.
\end{Proof}

We now have:
\begin{Proposition} If $\delta$ is not a positively marked Dynkin diagram, then $M(\delta)$ is not  orientable. \label{orientable}
\end{Proposition}
\begin{Proof} assume that the marked Dynkin diagram $\delta$ has some negative sign $s_{i,j}$. By Lemma \ref{twos} it is enough to compute any of the coefficients in (\ref{two3}) and show that  it is non-zero for the case of (\ref{gamma2}), the unstable manifold for $w=e$, the top cell.  By picking a negatively marked edge as close as possible to one of the ends of the Dynkin diagram we can reduce the calculation to the rank two cases.  This is an easy calculation which we illustrate in Example  (\ref{nonorientableA2}).
\end{Proof}

\begin{Example}
Non-orientability in the standard $A_2$ case.
\label{nonorientableA2}
\end{Example}
We compute one of the coefficients in  (\ref{two3}) . Consider $w=e$, $\Pi^s(w)=\emptyset$, $r=1$,  $\epsilon=(0,1)\in  {\cal E}^{\{\alpha_1 \}}$, $\epsilon(-1)=(-1,1)$ and $\epsilon(1)=(1,1)$. Now for $\sigma=s_{\alpha_1}$, we obtain $\sigma (\epsilon (-1))=(-1,-1)$ and $\sigma (\epsilon (1))=(1,1)$. Thus the product of all these signs is one and we have a non-zero boundary of the top cell (\ref{gamma2}).

\vskip 0.5cm

\begin{Corollary} The closure of an unstable manifold
  $\overline {C}^u(w)$ gives rise to a cycle in homology if and only if its marked Dynkin diagram $\delta_{\Pi^u(w)}$ is positively marked. \label{unstablecycle}
\end{Corollary}
\begin{Proof} We use Proposition (\ref{smoothunstable}) and  Proposition (\ref{orientable}) to conclude that  $\overline {C}^u(w)$ is orientable-thus a cycle- exactly when $\delta_{\Pi^u(w)}$ is positively marked.
\end{Proof}
\vskip 0.5cm

We now obtain in the standard case :
\begin{Corollary}If $M(\delta)$ corresponds to the standard marked Dynkin diagram, then  $\overline {C}^u(w)$ gives rise to a cycle in homology if and only if the unstable Weyl group $W^u(w)$ is abelian.
\end{Corollary}
\begin{Proof} In this case one observes that all the subdiagrams of standard marked Dynkin diagram remain standard. Hence these are not positively marked except in the case when there are no edges in the subdiagram and the corresponding  Weyl group $W^u(w)$ is abelian. 
\end{Proof}

Let $e(k)$ be the number of all elements in $W$ with index  $k$.
For example, $e(k)$ is the Eulerian number $E(l+1,k)$ for ${\frak g}$ of 
type $A_l$ (i.e. $W=S_{l+1}$, the symmetry group of order $l+1$).  Then we have the following theorem which is a generalization to the manifolds $M(\delta)$ 
considered in \cite{davis:91}.

\begin{Theorem} Let $\delta$ be a marked Dynkin diagram. Then for any  $k=0,1, \cdots ,l$  the homology group $H_k(M(\delta), {\Bbb Z}/2{\Bbb Z} )$ has rank $e(k)$.
\end{Theorem}
\begin{Proof} This follows from Morse theory and from Lemma (\ref {twos}) which implies that all the unstable manifolds are cycles over  ${\Bbb Z}/2{\Bbb Z}$.
\end{Proof}

Similarly, in the case of a Tomei manifold all the closures of unstable manifolds are  orientable since the marked Dynkin diagrams $\delta_{\Pi^u(w)}$ of Corollary \ref{unstablecycle} are positively marked.

We then recover a theorem of Davis \cite{davis:87}  (general setting of Coxeter groups) and Fried  \cite{fried:86} (for type $A_l$):

\begin{Theorem}  In the case when $M(\delta)$ is a Tomei manifold, for any  $k=0,1, \cdots ,l$ the homology group $H_k(M(\delta), {\Bbb Z} )$ is free of rank $e(k)$.
\end{Theorem}

\bibliographystyle{amsplain}

\end{document}